\documentstyle[12pt]{article}
\newcommand{\curl}{\mathop{\rm curl}\limits}
\newcommand{\const}{\mathop{\rm const}\limits}
\newcommand{\supp}{\mathop{\rm supp}\limits}
\newcommand{\mes}{\mathop{\rm mes}\limits}

\begin{document}

\begin{center}

{\bf Multidimensional Lusin-type inequalities }\\

\vspace{3mm}

  {\bf for Grand Lebesgue  Spaces. }\\

\vspace{3mm}

{\sc Ostrovsky E., Sirota L.}\\

\normalsize

\vspace{3mm}
{\it Department of Mathematics and Statistics, Bar-Ilan University,
59200, Ramat Gan, Israel.}\\
e-mail: \ eugostrovsky@list.ru \\

{\it Department of Mathematics and Statistics, Bar-Ilan University,
59200, Ramat Gan, Israel.}\\
e-mail: \ sirota3@bezeqint.net \\

\end{center}

\vspace{4mm}

  {\it Abstract.} We generalize in this short paper the classical  Luzin's theorem  about existence
of integral on the measurable function and its multidimensional analogues on the
many popular classes of rearrangement invariant (r.i.) spaces, namely, on the so-called Grand Lebesgue Spaces.\par

\vspace{4mm}

 \par

 \vspace{4mm}


\vspace{3mm}

\section{ Introduction. Notations. Statement of problem.}

\vspace{3mm}

 "{\it A classical theorem of Lusin \cite{Luzin1} states that for every Borel function $ f $ on $ R,$ there is a continuous function
$ u $ on $ R $ that is differentiable almost everywhere with derivative equal to $ f. $
In \cite{Alberti1}  G.Alberti gave a related result in higher dimensions. He proved the following theorem, in which $ | A | $
denotes the Lebesgue measure of the (Borelian) set $  A, \ A \subset R^d, \ d = 1,2,\ldots  $ and $ Du $ denotes the standard
derivative (gradient) of $ u. $ \par

\vspace{3mm}

 Theorem  (\cite{Alberti1}, Theorem 1). Let $ \Omega \subset R^k $ be open with $ |\Omega| < \infty, $ and let $ f : \Omega \to R^k $ be
a Borel function. Then for every $ \epsilon > 0, $ there exist an open set $ A \subset \Omega $ and a function $ u \in C^1 _0(\Omega) $
 such that

$$
 \ |A| \le \epsilon|\Omega|, \eqno(1.a)
$$

$$
 \ f(x) = Du(x), \ x \in  \Omega \setminus A, \eqno(1.b)
$$
and

$$
 \ |Du|_p \le K(d) \ \epsilon^{1/p - 1}|f|_p. \eqno(1.c)
$$

Here $ p \in (1, \infty), \  K(d) \in (0, \infty) $ is a constant that depends only on the dimension $ d, $

$$
|f|_p \stackrel{def}{=} \left[  \int_{\Omega}|f(x)|^p \ dx   \right]^{1/p}. \eqno(1.0)
$$

 In other words, Alberti showed that it is possible to arbitrarily prescribe the gradient of a $ C^1_0 $
function $ u $ on $ \Omega \in R^d $ off of a set of arbitrarily small measure, with quantitative control on all $ L_p $ norms of } $ Du", $
see \cite{David1}. \par

 Let us denote such a function $ u(\cdot), $ i.e. which satisfies the relation (1.c), not necessary to be unique, by

$$
u(x) = u_{\epsilon,A}(x) = u_{\epsilon,A}[f](x).
$$

 We can and will understood as a capacity of the "constant" $ K(d) $ its minimal  value

$$
K(d) := \sup_{p \in [1,\infty]}\sup_{\epsilon > 0} \sup_{ 0 \ne f \in L(p) }
 \left\{ \frac{|Du_{\epsilon,A}[f]|_p}{ \epsilon^{1/p - 1}|f|_p } \right\}. \eqno(1.1)
$$

 It follows immediately from the article \cite{Alberti1} after simple calculations that

$$
K(d) \le 72 \ d^{3/2}. \eqno(1.2)
$$

 Notice that the condition $  \curl(f) = 0,  $ where $ \curl(f) $  denotes the  distributional rotor of the function $  f,  $ is
not necessary for theorem 1! \par

 Note in addition that Alberti \cite{Alberti1} proved that the rate $ \epsilon^{1/p - 1} $ is sharp as $ \epsilon \to 0+.  $\par

See also the articles \cite{David1}, \cite{Cheeger1}, \cite{Moonens1}. In particular, the work  \cite{David1} contains
the extension of the results of Alberti \cite{Alberti1} and Moonens-Pfeffer \cite{Moonens1}, in a suitable sense, to a class
of metric measure spaces on which differentiation is defined. \par

\vspace{3mm}

{\bf  Our purpose in this short report is to show that the assertion of theorem 1 can be easily generalized
on the  wide class of the so - called moment rearrangement invariant spaces, namely, on the Grand Lebesgue Spaces (GLS).  } \par

\vspace{3mm}

 The so - called Grand Lebesgue Spaces (GLS) are very popular in recent years, see. e.g.
\cite{Capone1}, \cite{Fiorenza1}, \cite{Fiorenza2}, \cite{Iwaniec1},  \cite{Iwaniec2}, \cite{Jawerth1}, \cite{Kozatchenko1},
\cite{Liflyand1}, \cite{Ostrovsky1}, \cite{Ostrovsky2}. \\

  The  Grand Lebesgue Spaces $ GLS = G(\psi) =G\psi =
 G(\psi; A,B), \ A,B = \const, A \ge 1, A < B \le \infty, $ spaces consist by definition
 on all the measurable functions $ f: D \to R $ with finite norms

$$
   ||f||G(\psi) \stackrel{def}{=} \sup_{p \in (A,B)} \left[ |f|_p /\psi(p) \right].\eqno(1.3)
$$

  Here $ \psi(\cdot) $ is some continuous positive on the {\it open} interval
$ (A,B) $ function such that

$$
     \inf_{p \in (A,B)} \psi(p) > 0, \ \psi(p) = \infty, \ p \notin (A,B).
$$
 We will denote
$$
 \supp (\psi) \stackrel{def}{=} (A,B) = \{p: \psi(p) < \infty, \} \eqno(1.4)
$$

The set of all $ \psi $  functions with support $ \supp (\psi)= (A,B) $ will be
denoted by $ \Psi(A,B). $ \par
  This spaces are rearrangement invariant, see \cite{Bennet1}, and
    are used, for example, in the theory of probability  \cite{Kozatchenko1},
  \cite{Ostrovsky1}, \cite{Ostrovsky2}; theory of Partial Differential Equations
    \cite{Iwaniec2};  functional analysis  \cite{Iwaniec2},  \cite{Liflyand1},
  \cite{Ostrovsky2}; theory of Fourier series, theory of martingales, mathematical statistics,
  theory of approximation etc.\par

 Notice that in the case when $ \psi(\cdot) \in \Psi(A,\infty)  $ and a function
 $ p \to p \cdot \log \psi(p) $ is convex,  then the space
$ G\psi $ coincides with some {\it exponential} Orlicz space. \par
 Conversely, if $ B < \infty, $ then the space $ G\psi(A,B) $ does  not coincides with
 the classical rearrangement invariant spaces: Orlicz, Lorentz, Marcinkiewicz  etc.\par

  The fundamental function of these spaces $ \phi(G(\psi), \delta) = ||I_A ||G(\psi), \mes(A) = \delta, \ \delta > 0, $
where  $ I_A  $ denotes as ordinary the indicator function of the measurable set $ A, $ may be calculated by the formulae

$$
\phi(G(\psi), \delta) = \sup_{ p \in \supp (\psi)} \left[ \frac{\delta^{1/p}}{\psi(p)} \right]. \eqno(1.5)
$$

 The fundamental function of arbitrary rearrangement invariant spaces plays very important role in functional analysis,
theory of Fourier series and transform \cite{Bennet1} as well as in our further narration. \par

 Many examples of fundamental functions for some $ G\psi $ spaces are calculated and investigated in  \cite{Ostrovsky1}, \cite{Ostrovsky2}.\par

\vspace{3mm}

{\bf Remark 1.1} If we introduce the {\it discontinuous} function $ \psi_{(r)} =  \psi_{(r)}(p) $ such that

$$
\psi_{(r)}(p) = 1, \ p = r; \hspace{5mm} \psi_{(r)}(p) = \infty, \ p \ne r, \ p,r \in (A,B) \eqno(1.6)
$$
and define formally  $ C/\infty = 0, \ C = \const \in R^1, $ then  the norm
in the space $ G(\psi_r) $ coincides with the $ L_r $ norm:

$$
||f||G(\psi_{(r)}) = |f|_r. \eqno(1.7)
$$
Thus, the Grand Lebesgue Spaces are direct generalization of the
classical exponential Orlicz's spaces and Lebesgue - Riesz spaces $ L_r. $ \par

 Let a function $ f : \Omega \to R $ be such that

$$
\exists s_1, s_2, \ 1 \le s_1 < s_2 \le \infty: \ \forall p \in (s_1, s_2) \ \Rightarrow |f|_p < \infty.
$$

 Then the function $ \psi = \psi_f(p), \ s_1 < p < s_2  $ may be {\it naturally} defined by the following way:

$$
\psi_f(p)  := |f|_p; \hspace{6mm} \supp \psi_f(\cdot) = (s_1,s_2). \eqno(1.7)
$$

\vspace{3mm}

\section{ Main result.}

\vspace{3mm}

 Suppose $  \exists (A,B) = \const, 1  \le A < B \le \infty,   $ such that

 $$
\exists \psi(\cdot) \in G\Psi(A,B) \ \Rightarrow f(\cdot) \in G\psi. \eqno(2.1)
 $$
 For instance, let the (measurable) function $ f: \Omega \to R $ be such that

 $$
 \exists (A,B) = \const, \ 1 \le A < B \le \infty,  \ \forall p \in (A,B) \Rightarrow \ |f|_p < \infty. \eqno(2.2)
 $$
Then the function $  \psi = \psi(p) $ can be picked as a natural function for the function $  f: \  \psi(p) = \psi_f(p). $   \par

 Let also $ \zeta = \zeta(p)  $ be another $ \psi \ -  $ function with at the same support $ \supp \zeta = (A,B);  $ define the
new function $ \nu(p) = \zeta(p) \ \psi(p).  $  Obviously, $ \nu(\cdot) \in \Psi(A,B). $\par

\vspace{3mm}

{\bf Theorem 2.1.} {\sc We propose under notations and conditions of Theorem 1}

$$
||Du||G(\nu) \le K(d) \ \epsilon^{-1} \ \phi(G\zeta, \epsilon) \ ||f||G\psi, \eqno(2.3)
$$
{\sc where the defined before in (1.1) constant  $  K(d) $ in the inequality (2.3) is the best possible.} \par

\vspace{3mm}

{\bf Proof.} We can and will suppose without loss of generality $ ||f||G\psi = 1. $ The last relation implies in particular that

$$
\forall p \in (A,B) \ \Rightarrow |f|_p \le \psi(p). \eqno(2.4)
$$

 We derive substituting into (1.c)

$$
  |Du|_p \le K(d) \ \epsilon^{1/p - 1} \ \psi(p), \ p \in (A,B), \eqno(2.5)
$$
or equally

$$
\frac{\epsilon}{K(d)} \cdot \frac{|Du|_p}{\nu(p)} \le    \frac{\epsilon^{1/p}}{\zeta(p)}, \ p \in (A,B). \eqno(2.6)
$$

 We obtain taking supremum over $ p, \ p \in (A,B)  $
 from both the sides of inequality (2.6) taking into account the direct definition of the
Grand Lebesgue Space norm as well as the direct  definition of fundamental function of these spaces

$$
\frac{\epsilon}{K(d)} \cdot ||Du||G(\nu)  \le \phi(G\zeta, \epsilon) = \phi(G\zeta, \epsilon) \cdot ||f||G\psi,
$$
which is entirely equivalent to the first assertion of theorem 2.1. \par

 The sharpness  of the constant $  K(d) $ follows immediately from  one of the results of  the preprint
\cite{Ostrovsky3}. \par

 This completes the proof of our theorem.\par

\vspace{3mm}

{\bf Remark 2.1.} If we choose for instance  as a capacity of the function $  \zeta(p) $ the degenerate function
$  \zeta(p) = \psi_{(r)}(p),  \ r = \const > 1, $ and agree to take

$$
\psi_{(r)} \cdot \psi(p) = \psi_{(r)} \cdot \psi(r), \eqno(2.7)
$$
we obtain the Alberti's proposition (1.c) as a particular case. \par

\end{document}